\documentclass{amsart}
\usepackage{amssymb}
\newtheorem{thm}{Theorem}
\newtheorem{lem}{Lemma}
\newtheorem{prop}{Proposition}
\newtheorem{cor}{Corollary}

\newtheorem{df}{Definition}

\begin{document}

\bibliographystyle{plain}

\title[Milo\v s Arsenovi\'c and Romi F. Shamoyan]{On some new theorems on multipliers in harmonic function spaces in higher dimension II}

\author[]{Milo\v s Arsenovi\' c$\dagger$}
\author[]{Romi F. Shamoyan}

\address{Department of mathematics, University of Belgrade, Studentski Trg 16, 11000 Belgrade, Serbia}
\email{\rm arsenovic@matf.bg.ac.rs}

\address{Bryansk University, Bryansk Russia}
\email{\rm rshamoyan@yahoo.com}

\thanks{$\dagger$ Supported by Ministry of Science, Serbia, project OI174017}

\date{}

\begin{abstract}
We present new sharp assertions concerning multipliers in various spaces of harmonic functions in the unit ball of $\mathbb R^n$.
\end{abstract}

\maketitle

\footnotetext[1]{Mathematics Subject Classification 2010 Primary 30H20.  Key words
and Phrases: multipliers, spaces of harmonic functions, Bergman type mixed norm spaces, spherical harmonics.}

\section{Introduction and preliminaries}

In this paper we continue investigation, started in \cite{AS}, of spaces of multipliers between certain spaces of harmonic functions on the unit ball. While the subject of multipliers between spaces of analytic functions in the unit disc is a vast one, multipliers between spaces of analytic functions in the unit ball and in the unit polydisc in $\mathbb C^n$ are less explored, for some results in this direction see \cite{Sh1} and references therein. For results on multipliers between harmonic spaces on the unit disc we refer the reader to \cite{SW}, where multipliers between harmonic Bergman type classes were considered, and to \cite{Pa1} and \cite{Pa2} for the case of harmonic Hardy classes.

Let us briefly describe the content of the paper. In this section we describe spaces of harmonic functions on the
unit ball that are of interest to us and recall definition of multipliers between two such spaces. The next section
contains auxiliary results, however some embedding results presented there could be of independent interest. The last
section begins with general necessary conditions for a sequence to be a multiplier, these are valid for quite general
assumptions on the parameters involved in the definition of the spaces. In many cases these necessary conditions turn out to be sufficient as well and these characterizations of multipliers between various spaces of harmonic functions are the main results of this paper.

Let $\mathbb B$ be the open unit ball in $\mathbb R^n$, $\mathbb S = \partial \mathbb B$ is the unit sphere in $\mathbb R^n$, for $x \in \mathbb R^n$ we have $x = rx'$, where $r = |x| = \sqrt{\sum_{j=1}^n x_j^2}$ and $x' \in \mathbb S$. Let $\omega_n$ denote the volume of the unit ball in $\mathbb R^n$. Normalized Lebesgue measure on $\mathbb B$ is denoted by $dx = dx_1 \ldots dx_n = r^{n-1}dr dx'$ so that $\int_{\mathbb B} dx = 1$. We set $I = [0, 1)$. We denote the space of all harmonic functions in an open set $\Omega$ by $h(\Omega)$. The gradient of $f \in C^1(\Omega)$ is denoted by $\nabla f$, $|\nabla f(x)| = \sqrt{\sum_{j=1}^n |\partial f(x)/\partial x_j|^2}$ is its Euclidean norm. In this paper letter $C$ designates a positive constant which can change its value even in the same chain of inequalities.

For $0<p<\infty$, $0 \leq r < 1$ and $f \in h(\mathbb B)$ we set
$$M_p(f, r) = \left( \int_{\mathbb S} |f(rx')|^p dx' \right)^{1/p},$$
with the usual modification to cover the case $p = \infty$. For $0 < p < \infty$, $0 < q \leq \infty$, $\alpha > 0$
and $f \in h(\mathbb B)$ we consider mixed (quasi)-norms $\| f \|_{p,q,\alpha}$ defined by
\begin{equation}\label{qpnorm}
\| f \|_{p,q,\alpha} = \left( \int_0^1 M_q(f, r)^p (1-r^2)^{\alpha p - 1} r^{n-1}dr \right)^{1/p},
\end{equation}
and the corresponding spaces
$$B^{p,q}_\alpha(\mathbb B) = B^{p,q}_\alpha = \{ f \in h(\mathbb B) : \| f \|_{p,q,\alpha} < \infty \}.$$
For details on these spaces we refer to \cite{DS}, Chapter 7. In particular these spaces are complete metric spaces and for $\min(p,q) \geq 1$ they are Banach spaces. These mixed norm spaces include weighted Bergman spaces $A^p_\beta(\mathbb B) = A^p_\beta = B^{p,p}_{\frac{\beta +1}{p}}$ where $\beta > -1$ and $0 < p < \infty$, see \cite{DS} for more on these
spaces. We set $A^\infty_\beta = B^{\infty, \infty}_\beta$ for $\beta > 0$. We also consider, for $\alpha \geq 0$, weighted Hardy spaces
$$H^p_\alpha (\mathbb B) = H^p_\alpha = \{ f \in h(\mathbb B) : \| f \|_{p, \alpha} = \sup_{r<1} M_p(f, r)(1-r)^\alpha
< \infty \}, \qquad 0 < p \leq \infty.$$
If $\alpha = 0$ we write simply $H^p$ for $H^p_0$. In view of the above definitions it is natural to alow $p = \infty$
in the definition of $B^{p,q}_\alpha$ by setting $B^{\infty, p}_\alpha = H^p_\alpha$.

We denote by $B$ the harmonic Bloch space, i.e. the space of all functions $f \in h(\mathbb B)$ with finite norm
$$\| f \|_B = |f(0)| + \sup_{x \in \mathbb B} (1-|x|^2)|\nabla f(x)|.$$
This space is a Banach space, its closed subspace consisting of all $f \in B$ such that $\lim_{x \to \mathbb S}
(1-|x|^2)|\nabla f(x)| = 0$ is called little Bloch space and denoted by $B_0$. We refer the reader to \cite{JP2} for
details on $B$ and $B_0$.

Analytic Triebel-Lizorkin spaces were studied by several authors, see for example \cite{OF1}, \cite{Sh2} and references therein. In this paper we consider harmonic Triebel-Lizorkin spaces on the unit ball in $\mathbb R^n$.
\begin{df}
Let $0 < p,q < \infty$ and $\alpha > 0$. The harmonic Triebel-Lizorkin space $F^{p,q}_\alpha (\mathbb B) =
F^{p,q}_\alpha$ consists of all functions $f \in h(\mathbb B)$ such that
\begin{equation}\label{tl}
\| f \|_{F^{p,q}_\alpha} = \left( \int_{\mathbb S} \left( \int_0^1 |f(rx')|^p (1-r)^{\alpha p -1} dr \right)^{q/p}
dx' \right)^{1/q} < \infty.
\end{equation}
\end{df}
These spaces are complete metric spaces, for $\min(p,q) \geq 1 $ they are Banach spaces. We prove certain inclusions
between $F^{p,q}_\alpha$ and $B^{p,q}_\alpha$ spaces, see Propositions \ref{fembb} and \ref{byal}.

Next we need certain facts on spherical harmonics, see \cite{StW} for a detailed exposition. Let $Y^{(k)}_j$ be the spherical harmonics of order $k$, $j \leq 1 \leq d_k$, on $\mathbb S$. Next,
$$Z_{x'}^{(k)}(y') = \sum_{j=1}^{d_k} Y_j^{(k)}(x') \overline{Y_j^{(k)}(y')}$$
are zonal harmonics of order $k$. The spherical harmonics $Y^{(k)}_j$, ($k \geq 0$, $1 \leq j \leq d_k$), form an orthonormal basis of $L^2(\mathbb S, dx')$. Every $f \in h(\mathbb B)$ has an expansion
$$f(x) = f(rx') = \sum_{k=0}^\infty r^k b_k\cdot Y^k(x'),$$
where $b_k = (b_k^1, \ldots, b_k^{d_k})$, $Y^k = (Y_1^{(k)}, \ldots, Y_{d_k}^{(k)})$ and $b_k\cdot Y^k$ is interpreted in the scalar product sense: $b_k \cdot Y^k = \sum_{j=1}^{d_k} b_k^j Y_j^{(k)}$. We often write, to stress dependence on a function $f \in h(\mathbb B)$, $b_k = b_k(f)$ and $b_k^j = b_k^j(f)$, in fact we have linear functionals $b_k^j$,
$k \geq 0, 1 \leq j \leq d_k$, on the space $h(\mathbb B)$.

We denote the Poisson kernel for the unit ball by $P(x, y')$, it is given by
\begin{align*}
P(x, y') = P_{y'}(x) & = \sum_{k=0}^\infty r^k \sum_{j=1}^{d_k} Y^{(k)}_j(y') Y^{(k)}_j(x') \\
& = \frac{1}{n\omega_n} \frac{1-|x|^2}{|x-y'|^n}, \qquad x = rx' \in \mathbb B, \quad y' \in \mathbb S.
\end{align*}

We recall some definitions from \cite{AS}, these are needed to formulate our main results.

\begin{df}
For a double indexed sequence of complex numbers
$$c = \{ c_k^j : k \geq 0, 1 \leq j \leq d_k \}$$
and a harmonic function $f(rx') = \sum_{k=0}^\infty r^k \sum_{j=1}^{d_k} b_k^j(f) Y^{(k)}_j(x')$ we define
$$(c \ast f) (rx') = \sum_{k=0}^\infty \sum_{j=1}^{d_k} r^k c_k^j b_k^j(f) Y^{(k)}_j(x'), \qquad rx' \in \mathbb B,$$
if the series converges in $\mathbb B$. Similarly we define convolution of $f, g \in h(\mathbb B)$ by
$$(f \ast g)(rx') = \sum_{k=0}^\infty \sum_{j=1}^{d_k} r^k b_k^j(f)b_k^j(g) Y_j^{(k)}(x'), \qquad rx'\in \mathbb B,$$
it is easily seen that $f \ast g$ is defined and harmonic in $\mathbb B$.
\end{df}

\begin{df}
For $t > 0$ and a harmonic function $f(x) = \sum_{k=0}^\infty r^k b_k(f)\cdot Y^k(x')$ on $\mathbb B$ we define a
fractional derivative of order $t$ of $f$ by the following formula:
$$(\Lambda_t f)(x) = \sum_{k=0}^\infty r^k \frac{\Gamma(k+n/2 + t)}{\Gamma(k+n/2)\Gamma(t)}b_k(f)\cdot Y^k(x'),
\qquad x = rx' \in \mathbb B.$$
\end{df}
Clearly, for $f \in h(\mathbb B)$ and $t>0$ the function $\Lambda_t h$ is also harmonic in $\mathbb B$. We also note that $(g \ast P_{y'})(rx') = (g \ast P_{x'})(ry')$ and $\Lambda_t(f \ast g)(x) = (\Lambda_t f \ast g)(x)$ for any
$f, g \in h(\mathbb B)$, these easy to prove formulae are often used in our proofs without comment.

\begin{df}
Let $X$ and $Y$ be subspaces of $h(\mathbb B)$. We say that a double indexed sequence $c$ is a multiplier from $X$ to $Y$ if $c \ast f \in Y$ for every $f \in X$. The vector space of all multipliers from $X$ to $Y$ is denoted by $M_H(X, Y)$.
\end{df}

Clearly every multiplier $c \in M_H(X, Y)$ induces a linear map $M_c : X \rightarrow Y$. If, in addition, $X$ and $Y$ are complete (quasi)-normed spaces such that all functionals $b_k^j$ are continuous on both spaces $X$ and $Y$, then the map $M_c : X \rightarrow Y$ is continuous, as is easily seen using the Closed Graph Theorem. This condition is satisfied by all spaces we considered above.

\section{Auxiliary results}

In this section we collect results needed for our main results in the next section. In proving necessary conditions
for a double indexed sequence to be a multiplier one uses test functions, these are provided by the Bergman kernel
for harmonic weighted Bergman spaces $A^p_m$, $m>-1$. This kernel is the following function
\begin{equation}\label{bker}
Q_m(x, y) = 2 \sum_{k=0}^\infty \frac{\Gamma(m + 1 + k + n/2)}{\Gamma(m + 1) \Gamma(k + n/2)}
r^k \rho^k Z_{x'}^{(k)}(y'), \qquad x = rx', \; y = \rho y' \in \mathbb B.
\end{equation}
The test functions we are going to use in the next section are harmonic functions $f_{m, y}(x) = Q_m(x, y)$,
$y \in \mathbb B$. We often write $f_y$ instead of $f_{m, y}$.

The following lemma gives an estimate for the kernel $Q_m$, see \cite{DS}, \cite{JP1}.

\begin{lem}\label{DSlemma}
1. Let $m > -1$. Then, for $x = rx', y = \rho y' \in \mathbb B$ we have
$$|Q_m(x, y)| \leq \frac{C}{|\rho x - y'|^{n+m}}.$$
2. Let $m > n-1$, $0 \leq r < 1$ and $y'\in \mathbb S$. Then
$$\int_{\mathbb S} \frac{dx'}{|rx' - y'|^m} \leq \frac{C}{(1-r)^{m-n+1}}.$$
\end{lem}
The following lemma is an often used result from \cite{DS}.

\begin{lem}[\cite{DS}]\label{rro}
Let $\alpha > -1$ and $\lambda > \alpha + 1$. Then
$$\int_0^1 \frac{(1-r)^\alpha}{(1-r\rho)^\lambda} dr \leq C (1-\rho)^{\alpha + 1 - \lambda}, \qquad 0 \leq \rho < 1.$$
\end{lem}

We need some norm estimates of $f_y$.

\begin{lem}\label{estim}
Let $0 < t \leq \infty$. Then we have
$$M_t(f_{m,y}, r) \leq C (1-r|y|)^{-n-m + \frac{n-1}{t}}, \qquad m > \max \left(\frac{n-1}{t} - n, -1\right).$$
\end{lem}

This follows immediately from Lemma \ref{DSlemma}. As a consequence we obtain, using Lemma \ref{rro}, the following proposition.

\begin{prop}\label{prbpq}
Let $0 < p < \infty$, $0 < t \leq \infty$. Then, for $m > \max (\alpha + \frac{n-1}{t} - n, -1)$ we have:
$$\| f_{m,y} \|_{B^{p,t}_\alpha} \leq C (1-|y|)^{\alpha - n - m + \frac{n-1}{t}}, \qquad y \in \mathbb B.$$
\end{prop}

The missing case $p = \infty$, i.e. the case of Hardy spaces, is treated in the next proposition.
\begin{prop}
Let $0 < t \leq \infty$, $\alpha \geq 0$ and $m > \max (\alpha -n +\frac{n-1}{t}, -1)$. Then we have
\begin{equation}\label{hqnorm}
\| f_{m,y} \|_{H^t_\alpha} \leq C(1-|y|)^{\alpha -n-m+\frac{n-1}{t}}, \qquad y \in \mathbb B.
\end{equation}
\end{prop}

The proof of this proposition is similar to the proof of estimate (11) from \cite{AS} and can be left to the reader.
The following lemma is a preparation for analogous estimates of $F^{p,t}_\alpha$ norm of $f_{m,y}$.

\begin{lem}
Let $\gamma > -1$, $0<p<\infty$, $m > -1$ and $p(n+m) > \gamma + 1$. Then we have
\begin{equation}
\int_0^1 |f_{m,y}(rx')|^p (1-r)^\gamma dr \leq C |x' - y|^{\gamma + 1 - p(n+m)}, \qquad y \in \mathbb B, \quad x'
\in \mathbb S.
\end{equation}
\end{lem}

{\it Proof.} Using definitions and Lemma \ref{DSlemma} we obtain
\begin{equation*}
M  = \int_0^1 |f_{m,y}(rx')|^p (1-r)^\gamma dr \leq C \int_0^1 \frac{(1-r)^\gamma dr}{|x' - ry|^{p(n+m)}}.
\end{equation*}
Using elementary geometric inequality $|x' - ry| \geq c (|x' - y | + (1-r))$ we obtain, using Lemma \ref{rro},
$$M \leq C \int_0^1 \frac{(1-r)^\gamma dr}{[|x' - y| + (1-r)]^{p(n+m)}} \leq C |x' - y|^{\gamma + 1 -p(n+m)}.
\qquad \Box$$

\begin{prop}
Let $0<p,t<\infty$ and $m > \max (\alpha + \frac{n-1}{t} - n, -1)$. Then we have:
\begin{equation}\label{fnorm}
\| f_{m,y} \|_{F^{p,t}_\alpha} \leq C (1-|y|)^{\alpha - n - m + \frac{n-1}{t}}, \qquad y \in \mathbb B.
\end{equation}
\end{prop}

{\it Proof.} Using the above lemma and Lemma \ref{DSlemma} we obtain
\begin{align*}
\| f_{m,y} \|_{F^{p,t}_\alpha}^t & = \int_{\mathbb S} \left( \int_0^1 |f(rx')|^p (1-r)^{\alpha p -1} dr \right)^{t/p}
dx' \\
& \leq C \int_{\mathbb S} |x'-y|^{t(\alpha - n - m)} dx' \leq C (1-|y|)^{t(\alpha - n-m) + n-1}. \qquad \Box
\end{align*}

The remaining part of this section is devoted to embedding results, these are often used in proofs of our main
results in the next section.

\begin{lem}\label{qlo}
If $0<s \leq t \leq \infty$ then
\begin{equation}\label{eqglo}
M_t(f, r) \leq C(1-r)^{(n-1)(1/t-1/s)} M_s(f, r), \qquad 0 \leq r < 1, \quad f \in h(\mathbb B).
\end{equation}
\end{lem}

{\it Proof.} We can assume $s<t$. Let us set $b = (n-1)(\frac{1}{s} - \frac{1}{t})$ and let $I_b$ denote the operator of fractional integration of order $b > 0$. Now (\ref{eqglo}) follows immediately from the following estimates:
\begin{align}
M_s(I_bf, r) & \leq C(1-r)^{-b} M_s(f,r), \qquad 0\leq r < 1, \\
M_t(f,r) & \leq CM_s(I_bf, r), \qquad 0 \leq r < 1.
\end{align}
For the first estimate see \cite{DS}, Chapter 7, for the second, which is a Hardy-Littlewood theorem for harmonic functions in the unit ball, see \cite{Al1}, page 47. $\Box$

If we raise both sides of inequality (\ref{eqglo}) to the $p$-th power, multiply by $(1-r)^{\alpha p -1}$ and integrate
over $I$ we obtain the following corollary.
\begin{cor}
Let $0<s \leq t<\infty$, $0<p<\infty$ and $\alpha > (n-1)(\frac{1}{s}-\frac{1}{t})$. Then
\begin{equation}\label{qloemb}
\| f \|_{B^{p,t}_\alpha} \leq C \| f \|_{B^{p,s}_\beta}, \qquad \beta = \alpha + (n-1)(1/t-1/s) \quad f \in
h(\mathbb B),
\end{equation}
i.e. $B^{p,s}_\beta$ is continuously embedded into $B^{p,t}_\alpha$.
\end{cor}

The case $p = \infty$ leads us to weighted Hardy spaces, for $0<s \leq t<\infty$ the following result is an immediate consequence of Lemma \ref{qlo}:
\begin{equation}\label{harglo}
\| f \|_{H^t_\alpha} \leq C \| f \|_{H^s_\beta}, \qquad \beta = \alpha + (n-1)(1/t-1/s) \quad f \in h(\mathbb B).
\end{equation}

Let us note another continuous embedding from \cite{DS}:
\begin{equation}\label{inc}
B^{p_0, s}_\alpha \hookrightarrow B^{p_1, s}_\alpha, \qquad 0<s\leq\infty, \quad 0< p_0 \leq p_1 \leq \infty, \quad
\alpha > 0.
\end{equation}

\begin{prop}\label{fembb}
For $0<t \leq p<\infty$ the space $F^{p,t}_\alpha$ is continuously embedded into $B^{p,t}_\alpha$:
\begin{equation}\label{eqfb}
\| f \|_{B^{p,t}_\alpha} \leq C \| f \|_{F^{p,t}_\alpha}, \qquad  \alpha > 0.
\end{equation}
\end{prop}

{\it Proof.} For $f \in B^{p,t}_\alpha$ we have, using continuous form of Minkowski's inequality
\begin{align*}
\| f \|_{B^{p,t}_\alpha}^p &
\leq C \int_0^1 \left( \int_{\mathbb S} |f(rx')|^t dx' \right)^{p/t} (1-r)^{\alpha p-1} dr\\
& = C \left\| (1-r)^{\alpha t} \int_{\mathbb S} |f(rx')|^t dx' \right\|^{p/t}_{L^{p/t}((1-r)^{-1}dr)} \\
& \leq C \left( \int_{\mathbb S} \left\|(1-r)^{\alpha t}|f(rx')|^t \right\|_{L^{p/t}((1-r)^{-1}dr)} dx' \right)^{p/t}\\
& = C \| f \|_{F^{p,t}_\alpha}^p. \qquad \Box
\end{align*}

\begin{prop}\label{byal}
Let $0<p\leq 1$ and $p \leq q$. Then $B^{p,p}_\alpha$ is continuously embedded into $F^{q,p}_\alpha$:
$$\| f \|_{F^{q,p}_\alpha} \leq C \| f \|_{B^{p,p}_\alpha}, \qquad f \in h(\mathbb B).$$
\end{prop}

{\it Proof.} We set $I_n = [1-2^{-n}, 1-2^{-n-1})$, then we have $\cup_{n=0}^\infty I_n = I$. Set, for
$x' \in \mathbb S$, $I_n(x') = \int_{I_n} |f(rx')|^q(1-r)^{\alpha q -1} dr$. Since $p/q \leq 1$ we have
\begin{align*}
\| f \|_{F^{q,p}_\alpha}^q & = \int_{\mathbb S} \left( \int_0^1 |f(rx')|^q (1-r)^{\alpha q-1} dr\right)^{p/q} dx'\\
& = \int_{\mathbb S} \left( \sum_{n=0}^\infty I_n(x') \right)^{p/q} dx' \leq \int_{\mathbb S} \sum_{n=0}^\infty
I_n(x')^{p/q} dx'\\
& = \sum_{n=0}^\infty \int_{\mathbb S} I_n(x')^{p/q} dx'.
\end{align*}
Set $M_nf(x') = \sup_{r \in I_n} |f(rx')|$ for $n \geq 0$ and $x' \in \mathbb S$. Then we have, see \cite{Al1},
page 47,
$$\int_{\mathbb S} M^p_nf(x') dx' \leq C 2^n \int_{I_n} M_p^p(f, r) dr.$$
Clearly, $I_n(x') \leq C 2^{-n\alpha q} M_n^qf(x')$ for $x' \in \mathbb S$ and therefore we can use the above estimate to get
\begin{align*}
\| f \|_{F^{q,p}_\alpha}^q & \leq C \sum_{n=0}^\infty \int_{\mathbb S} 2^{-np\alpha} M_n^pf(x') dx' \leq C \sum_{n=0}^\infty 2^{-n(p\alpha-1)} \int_{I_n} M_p^p(f, r) dr \\
& \leq C \sum_{n=0}^\infty \int_{I_n} M_p^p(f, r) (1-r)^{\alpha p -1} dr \\
& = \| f \|_{B^{p,p}_\alpha}. \qquad \Box
\end{align*}

\section{Multipliers on spaces of harmonic functions}

In this section we present our main results: sufficient and/or necessary conditions for a double indexed sequence $c$ to be in $M_H(X, Y)$, for certain (quasi) normed spaces $X$ and $Y$ of harmonic functions. We associate to such a sequence $c$ a harmonic function
\begin{equation}\label{gc}
g_c(x) = g(x) = \sum_{k\geq 0} r^k \sum_{j=1}^{d_k} c_k^j Y^{(k)}_j(x'), \qquad x = rx' \in \mathbb B,
\end{equation}
and express our conditions in terms of fractional derivatives of $g_c$.

The first part of the lemma below appeared, in dimension two, in \cite{Pa2}, for the second part see \cite{AS} and
\cite{SA}.

\begin{lem}\label{intcon}
Let $f, g \in h(\mathbb B)$ have expansions
\begin{equation*}
f(rx')  = \sum_{k=0}^\infty r^k \sum_{j=1}^{d_k} c_k^j Y^{(k)}_j(x'), \qquad
g(rx')  = \sum_{l=0}^\infty r^l \sum_{i=1}^{d_k} b_l^i Y^{(l)}_i(x').
\end{equation*}
Then we have
\begin{equation*}
\int_{\mathbb S} (g \ast P_{y'})(rx') f(\rho x') dx' = \sum_{k=0}^\infty r^k\rho^k \sum_{j=1}^{d_k} b_k^j c_k^j
Y^{(k)}_j(y'), \qquad y' \in \mathbb S, \quad 0 \leq r, \rho  < 1.
\end{equation*}
Moreover, for every $m > -1$, $y' \in \mathbb S$ and $0 \leq r, \rho < 1$ we have
\begin{equation*}
\int_{\mathbb S} (g \ast P_{y'})(rx') f(\rho x') dx' = 2 \int_0^1 \int_{\mathbb S} \Lambda_{m+1}(g \ast P_{y'})(rRx')
f(\rho Rx') (1-R^2)^m R^{n-1} dx' dR.
\end{equation*}
\end{lem}

We note for future use the following formula, contained in Lemma \ref{intcon}:
\begin{equation}\label{simple}
(c \ast f)(r^2x') = \int_{\mathbb S} (g_c \ast P_{y'})(rx') f(ry') dy', \qquad r \in I, \quad x' \in \mathbb B.
\end{equation}
Also, if $h_y = h_{m,y} = M_c f_{m,y}$ where $m>-1$ and $y = \rho y' \in \mathbb B$, then
\begin{equation*}
h_y(x) = \sum_{k\geq 0} r^k \rho^k \sum_{j=1}^{d_k} \frac{\Gamma(k+n/2+m+1)}{\Gamma(k+n/2)\Gamma(m+1)} c_k^j
Y^{(k)}_j(y') Y^{(k)}_j(x'), \quad x = rx' \in \mathbb B.
\end{equation*}
This gives the following formula which will be in constant use:
\begin{equation}\label{hay}
h_y(x) = \Lambda_{m+1} (g_c \ast P_{y'})(\rho x) \qquad y = \rho y' \in \mathbb B, \quad x \in \mathbb B.
\end{equation}
The first part of the following lemma, which gives necessary conditions for $c$ to be a multiplier, is based on \cite{AS}.

\begin{lem}\label{nec}
Let $0 < p,q,t \leq \infty$, $1 \leq s \leq \infty$ and $m > \max (\alpha + \frac{n-1}{t} - n, -1)$. Assume a double indexed sequence $c = \{ c_k^j : k \geq 0, 1 \leq j \leq d_k \}$ is a multiplier from $B^{p,t}_\alpha$ to $B^{q,s}_\beta$ and $g = g_c$ is defined in {\rm(\ref{gc})}. Then the following condition is satisfied:
\begin{equation}\label{mult}
T_s(g) = \sup_{0\leq\rho < 1} \sup_{y' \in \mathbb S} (1-\rho)^{m-\alpha + \beta+n-\frac{n-1}{t}}
\left( \int_{\mathbb S} |\Lambda_{m+1}(g \ast P_{x'})(\rho y')|^s dx'\right)^{1/s} < \infty,
\end{equation}
where the case $s = \infty$ requires usual modification.

Also, let $0 < p, t \leq \infty$, $0 < s \leq \infty$, $\alpha > 0$, $\beta \geq 0$ and
$m > \max (\alpha + \frac{n-1}{t} - n, -1)$. If a double indexed sequence $c = \{ c_k^j : k \geq 0, 1 \leq j \leq d_k \}$ is a multiplier from $B^{p,t}_\alpha$ to $H^s_\beta$, then the above function $g$ satisfies condition \rm{(\ref{mult})}.
\end{lem}

{\it Proof.}  Let $c \in M_H(B^{p,t}_\alpha, B^{q,s}_\beta)$, and assume both $p$ and $q$ are finite, the infinite cases require only small modifications. We have $\| M_c f\|_{B^{q,s}_\beta} \leq C \| f \|_{B^{p,t}_\alpha}$ for  $f$ in $B^{p,t}_\alpha$. Set $h_y = M_c f_y$, then we have
\begin{equation}\label{cont}
\| h_y \|_{B^{q,s}_\beta} \leq C \| f_y \|_{B^{p,t}_\alpha}.
\end{equation}
This estimate and Proposition \ref{prbpq} give
\begin{equation}\label{hynorm}
\| h_y \|_{B^{q,s}_\beta} \leq C (1-|y|)^{\alpha - m - n + \frac{n-1}{t}}, \qquad y \in \mathbb B.
\end{equation}
Using (\ref{hay}) and monotonicity of $M_s(h_y, r)$ we obtain, for $y = \rho y' \in \mathbb B$:
\begin{align}\label{muka}
M_s(\Lambda_{m+1}(g \ast P_{y'}), \rho^2) & = \left( \int_{\rho}^1 (1-r)^{\beta q -1} r^{n-1} dr \right)^{-1/q}
\notag \\
& \phantom{=} \times \left( \int_\rho^1 (1-r)^{\beta q -1} r^{n-1} M_s^q(h_y, \rho^2) dr \right)^{1/q} \notag\\
& \leq C(1-\rho)^{-\beta} \left( \int_\rho^1 (1-r)^{\beta q -1} r^{n-1} M_s^q(h_y, r) dr \right)^{1/q} \notag\\
& \leq C(1-\rho)^{-\beta} \| h_y \|_{B^{q,s}_\beta}.
\end{align}
Combining (\ref{muka}) and (\ref{hynorm}) we obtain
\begin{equation*}
\left( \int_{\mathbb S} |\Lambda_{m+1}(g \ast P_{x'})(\rho^2 y')|^s dx'\right)^{1/s} \leq C
(1 - \rho)^{\alpha - \beta - m - n + \frac{n-1}{t}},
\end{equation*}
which is equivalent to (\ref{mult}). The case $s = \infty$ is treated similarly.

Next we consider $c \in M_H(B^{p,t}_\alpha, H^s_\beta)$, assuming $0<p \leq \infty$.
Set $h_y = M_c f_y = g \ast f_y$. We have, by Proposition \ref{prbpq},
\begin{equation}\label{is1}
\| f_y \|_{B^{p,t}_\alpha} \leq C(1-|y|)^{\alpha - m - n + \frac{n-1}{t}}, \qquad y \in \mathbb B,
\end{equation}
and, by continuity of $M_c$, $\| h_y \|_{H^s_\beta} \leq C \| f_y \|_{B^{p,t}_\alpha}$. Therefore
\begin{equation}\label{is2}
\| h_y \|_{H^s_\beta} \leq C(1-|y|)^{\alpha - m - n + \frac{n-1}{t}}, \qquad y \in \mathbb B.
\end{equation}
Setting $y = \rho y'$ we have
\begin{align}\label{muka1}
I_{y'}(\rho^2) & = \left( \int_{\mathbb S} |\Lambda_{m+1} (g \ast P_{x'})(\rho^2 y')|^s dx' \right)^{1/s}
= \left( \int_{\mathbb S} |\Lambda_{m+1} (g \ast P_y)(\rho x')|^s dx' \right)^{1/s} \notag \\
& = M_s(h_y, \rho) \leq (1-|y|)^{-\beta} \| h_y \|_{H^s_\beta}.
\end{align}
The last two estimates yield
$$\left( \int_{\mathbb S} |\Lambda_{m+1} (g \ast P_{x'})(\rho^2 y')|^s dx' \right)^{1/s} \leq C (1-|y|)^{\alpha - \beta
- m - n + \frac{n-1}{t}},\qquad |y| = \rho$$
which is equivalent to (\ref{mult}). $\Box$

The first part of the above lemma combined with Proposition \ref{fembb} gives the following corollary.
\begin{cor}\label{forbf}
If $c \in M_H(B^{p,t}_\alpha, F^{q,s}_\beta)$, where $0<p,t\leq \infty$, $1 \leq s \leq q < \infty$ and
$m > \max (\alpha + \frac{n-1}{t} - n, -1)$, then the function $g_c$ satisfies condition
\begin{equation*}
T_s(g) = \sup_{0\leq\rho < 1} \sup_{y' \in \mathbb S} (1-\rho)^{\beta - \alpha + m + n - \frac{n-1}{t}}
\left( \int_{\mathbb S} |\Lambda_{m+1}(g \ast P_{x'})(\rho y')|^s dx'\right)^{1/s} < \infty.
\end{equation*}
\end{cor}

\begin{lem}\label{nec1}
Let $0 < p, t < \infty$, $0 < q \leq \infty$, $1 \leq s \leq \infty$, and $m > \max (\alpha + \frac{n-1}{t} -n, -1)$. Assume $c \in M_H(F^{p,t}_\alpha, B^{q,s}_\beta)$. Then the function $g = g_c$ satisfies the following condition:
\begin{equation}\label{mult1}
T_s(g) = \sup_{0\leq\rho < 1} \sup_{y' \in \mathbb S} (1-\rho)^{m+\beta-\alpha+n-\frac{n-1}{t}} \left( \int_{\mathbb S}
|\Lambda_{m+1}(g \ast P_{x'})(\rho y')|^s dx'\right)^{1/s} < \infty.
\end{equation}

Next, if $c \in M_H(F^{p,t}_\alpha, H^s_\beta)$, where $0< p,t < \infty$, $0<s\leq\infty$, $\alpha > 0$, $\beta \geq 0$
and $m > \max (\alpha + \frac{n-1}{t} -n, -1)$, then $g = g_c$ satisfies condition {\rm (\ref{mult1})}.

Finally, if $c \in M_H(H^t_\alpha, H^s_\beta)$, where $m > \max (\alpha + \frac{n-1}{t} -n, -1)$, $0<s,t\leq \infty$ and
$\alpha \geq 0$, $\beta \geq 0$, then again $g = g_c$ satisfies condition {\rm (\ref{mult1})}.
\end{lem}

Proofs of all three statements are analogous to the proof of Lemma \ref{nec}. For the first one use (\ref{fnorm})
and (\ref{muka}), for the second one use (\ref{fnorm}) and (\ref{muka1}) and for the last one use (\ref{hqnorm}) and
(\ref{muka1}). We leave details to the reader.

In our previous work all characterizations of multipliers were independent on the dimension of the space. However, in the following theorem dimension of the space is present in the description of the space of multipliers.

\begin{thm}\label{btoh}
Let $0 < t,p \leq 1$, $1 \leq s \leq \infty$ and $m > \max (\alpha + \frac{n-1}{t} - n, -1)$. Then for a double indexed sequence $c = \{ c_k^j : k \geq 0, 1 \leq j \leq d_k \}$ the following conditions are equivalent:

1. $c \in M_H(B^{p,t}_\alpha, H^s_\beta)$.

2. The function $g(x) = \sum_{k\geq 0} r^k \sum_{j=1}^{d_k} c_k^j Y^{(k)}_j(x')$ is harmonic in $\mathbb B$ and
satisfies the following condition
\begin{equation}\label{ngdim}
T_s(g) = \sup_{0\leq\rho < 1} \sup_{y' \in \mathbb S} (1-\rho)^{\beta - \alpha + m +n-\frac{n-1}{t}}
\left( \int_{\mathbb S} |\Lambda_{m+1}(g \ast P_{x'})(\rho y')|^s dx'\right)^{1/s} < \infty.
\end{equation}
For $1 < t \leq \infty$ condition {\rm 1.} implies condition {\rm 2.}
\end{thm}

{\it Proof.} The necessity of the condition (\ref{ngdim}) is contained in Lemma \ref{nec}. Now we prove sufficiency of condition (\ref{ngdim}). Let $f \in B^{p,t}_\alpha$ and set
$h = M_c f$. We have, by Lemma \ref{intcon}:
\begin{equation*}\label{tool}
h(r^2 x') = 2 \int_0^1 \int_{\mathbb S} \Lambda_{m+1} (g \ast P_{\xi})(rR x') f(r R \xi) (1-R^2)^m R^{n-1} d\xi dR.
\end{equation*}
Therefore, since $s \geq 1$, we deduce
\begin{equation*}
M_s(h, r^2) \leq C T_s(g) \int_0^1 (1-R)^m M_1(f, rR) (1-rR)^{\alpha - \beta -m-n+\frac{n-1}{t}}R^{n-1} dR.
\end{equation*}
Now we use Lemma 3 from \cite{AS} and obtain
\begin{align*}
M_s^p(h, r^2) & \leq CT_s^p(g) \int_0^1
\frac{(1-R)^{mp + p -1}}{(1-rR)^{p(m+n+\beta-\alpha) - \frac{p}{t}(n-1)}}M_1^p(f, rR)R^{n-1} dR\\
& \leq CT_s^p(g) (1-r)^{-\beta p} \int_0^1
\frac{(1-R)^{mp + p -1}}{(1-rR)^{p(m+n-\alpha) - p\frac{n-1}{t}}}M_1^p(f, rR)R^{n-1} dR.
\end{align*}
This inequality, monotonicity of $M_1(f, r)$ and Lemma \ref{qlo} give
\begin{align*}
(1-r)^{\beta p} M_s^p(h, r^2) & \leq C T_s^p(g) \int_0^1 \frac{(1-R)^{mp + p-1}}{(1-rR)^{p(m+n-\alpha)-
p\frac{n-1}{t}}} M_1^p(f, R) R^{n-1}dR \\
& \leq CT_s^p(g) \int_0^1 (1-R)^{p-1+p(\frac{n-1}{t} + \alpha - n)} M_1^p(f,R) dR\\
& \leq CT_s^p(g) \int_0^1 (1-R)^{\alpha p-1} M_t^p(f, R)R^{n-1} dR\\
& = CT_s^p(g) \| f \|_{B^{p,t}_\alpha}^p,
\end{align*}
which implies $\| h \|_{H^s_\beta} \leq CT_s(g) \| f \|_{B^{p,t}_\alpha}$ and the proof is complete. $\Box$

Theorem \ref{btoh}, Lemma \ref{nec1} and Proposition \ref{fembb} combine to give the following result.
\begin{thm}
Let $0<t \leq p \leq 1 \leq s \leq \infty$ and $m > \max (\alpha + \frac{n-1}{t} - n, -1)$. Then for a double indexed sequence $c = \{ c_k^j : k \geq 0, 1 \leq j \leq d_k \}$ the following conditions are equivalent:

1. $c \in M_H(F^{p,t}_\alpha, H^s_\beta)$.

2. The function $g(x) = \sum_{k\geq 0} r^k \sum_{j=1}^{d_k} c_k^j Y^{(k)}_j(x')$ is harmonic in $\mathbb B$ and
satisfies the following condition
\begin{equation*}
T_s(g) = \sup_{0\leq\rho < 1} \sup_{y' \in \mathbb S} (1-\rho)^{\beta - \alpha + m +n-\frac{n-1}{t}}
\left( \int_{\mathbb S} |\Lambda_{m+1}(g \ast P_{x'})(\rho y')|^s dx'\right)^{1/s} < \infty.
\end{equation*}
\end{thm}

A characterization of $M_H(H^1_\alpha, H^s_\beta)$, $1 \leq s \leq \infty$, was given in \cite{AS}, the next two theorems generalize that result.
\begin{thm}\label{part1}
Let $0<t\leq 1 \leq s \leq \infty$, $\alpha \geq 0$, $m > \max (\alpha + \frac{n-1}{t} - n, -1)$ and $\beta > 0$. Then for a double indexed sequence $c = \{ c_k^j : k \geq 0, 1 \leq j \leq d_k \}$ the following conditions are equivalent:

1. $c \in M_H(H^t_\alpha, H^s_\beta)$.

2. The function $g(x) = \sum_{k\geq 0} r^k \sum_{j=1}^{d_k} c_k^j Y^{(k)}_j(x')$ is harmonic in $\mathbb B$ and
satisfies the following condition
\begin{equation}\label{ngh}
T_s(g) = \sup_{0\leq\rho < 1} \sup_{y' \in \mathbb S} (1-\rho)^{\beta - \alpha + m +n-\frac{n-1}{t}}
\left( \int_{\mathbb S} |\Lambda_{m+1}(g \ast P_{x'})(\rho y')|^s dx'\right)^{1/s} < \infty.
\end{equation}
\end{thm}

{\it Proof.} The necessity of (\ref{ngh}) is contained in Lemma \ref{nec1}. Let us prove sufficiency.
We choose $f \in H^t_\alpha$ and set $h = M_cf$. Applying the operator $\Lambda_{m+1}$ to equation (\ref{simple})
we obtain $\Lambda_{m+1} h(rx) = \int_{\mathbb S} \Lambda_{m+1}(g\ast P_{y'})(x) f(ry') dy'$. Since $s\geq 1$ this
gives:
\begin{align*}
M_s(\Lambda_{m+1}h, r^2) & \leq M_1(f, r) \sup_{y' \in \mathbb S} \| \Lambda_{m+1} (g \ast P_{y'})(x) \|_{L^s(dx')}\\
& \leq T_s(g) (1-r)^{\alpha - \beta -m-n+\frac{n-1}{t}} M_1(f,r) \\
& \leq T_s(g) (1-r)^{\alpha - \beta - m - 1} M_t(f,r),
\end{align*}
where we at the last step used Lemma \ref{qlo}. Hence we obtained
$$(1-r)^{m+1+\beta} M_s(\Lambda_{m+1} h, r^2) \leq C (1-r)^\alpha M_t(f, r), \qquad 0 \leq r < 1,$$
and, since $\beta > 0$, this implies $(1-r)^\beta M_s(h, r^2) \leq C \| f \|_{H^t_\alpha}$, see \cite{DS}, Chapter 7. Hence $\| h \|_{H^s_\beta} \leq C \| f \|_{H^t_\alpha}$. $\Box$

We note that the above proof of sufficiency does not work in the case $\beta = 0$. In the following theorem we deal with
unweighted Hardy spaces.
\begin{thm}\label{part1}
Let $0< t < 1 \leq s \leq \infty$ and $m > \max (\frac{n-1}{t} - n, -1)$. Then for a double indexed sequence $c = \{ c_k^j : k \geq 0, 1 \leq j \leq d_k \}$ the following conditions are equivalent:

1. $c \in M_H(H^t, H^s)$.

2. The function $g(x) = \sum_{k\geq 0} r^k \sum_{j=1}^{d_k} c_k^j Y^{(k)}_j(x')$ is harmonic in $\mathbb B$ and
satisfies the following condition
\begin{equation}\label{ngh1}
T_s(g) = \sup_{0\leq\rho < 1} \sup_{y' \in \mathbb S} (1-\rho)^{m +n-\frac{n-1}{t}}
\left( \int_{\mathbb S} |\Lambda_{m+1}(g \ast P_{x'})(\rho y')|^s dx'\right)^{1/s} < \infty.
\end{equation}
\end{thm}

{\it Proof.} As in the previous theorem, necessity of condition (\ref{ngh1}) follows from Lemma \ref{nec1}. Let $g = g_c$ satisfy (\ref{ngh1}), let $f \in H^t_\alpha$ and set $h = M_c f$. Then, using Lemma \ref{intcon} and continuous form of Minkowski's inequality we obtain:
\begin{align*}
M_s(h, r^2) & \leq C \int_0^1  M_s(\Lambda_{m+1}(g \ast P_{y'}), rR)\int_{\mathbb S} |f(rRx')| dx' (1-R)^m R^{n-1}dR \\
&\leq C T_s(g) \int_0^1 (1-rR)^{-m-n+\frac{n-1}{t}} (1-R)^m M_1(f, rR) dR\\
&\leq C T_s(g) \int_0^1 (1-R)^{-n+\frac{n-1}{t}} M_1(f, R) dR\\
&\leq C T_s(g) \| f \|_{H^t},
\end{align*}
the last estimate is a corollary of Carleson-Duren embedding theorem, see \cite{Al2}. $\Box$

The theorem below is the first result on multipliers into Triebel-Lizorkin spaces.
\begin{thm}
Let $0<p\leq 1 \leq q \leq \infty$ and $m > \alpha -1$. Then for a double indexed sequence
$c = \{ c_k^j : k \geq 0, 1 \leq j \leq d_k \}$ the following conditions are equivalent:

1. $c \in M_H(B^{p,1}_\alpha, F^{q,1}_\beta)$.

2. The function $g(x) = \sum_{k\geq 0} r^k \sum_{j=1}^{d_k} c_k^j Y^{(k)}_j(x')$ is harmonic in $\mathbb B$ and
satisfies the following condition
\begin{equation}\label{mult2}
N_1(g) = \sup_{0\leq\rho < 1} \sup_{y' \in \mathbb S} (1-\rho)^{\beta - \alpha + m + 1}
\int_{\mathbb S} |\Lambda_{m+1}(g \ast P_{x'})(\rho y')| dx' < \infty.
\end{equation}
\end{thm}

{\it Proof.} Necessity of condition (\ref{mult2}) is contained in Corollary \ref{forbf}. Now we choose $c$ such that
the condition (\ref{mult2}) is satisfied. Then, by Theorem 6 from \cite{AS} we have $M_c : B^{p,1}_\alpha \rightarrow B^{1,1}_\beta$. Since, by Proposition \ref{byal}, $B^{1,1}_\beta \hookrightarrow F^{q,1}_\beta$ the proof is
complete. $\Box$

The following theorem is a generalization of Theorem 3 from \cite{AS}. The proof is included for reader's convenience,
it follows almost the same pattern as in \cite{AS}.

\begin{thm}\label{bpbp}
Let $1 \leq p \leq q \leq \infty$, $1 \leq s \leq \infty$ and $m > \alpha - 1$. Then for a double indexed sequence
$c = \{ c_k^j : k \geq 0, 1 \leq j \leq d_k \}$ the following conditions are equivalent:

1. $c \in M_H(B^{p,1}_\alpha, B^{q,s}_\beta)$.

2. The function $g(x) = \sum_{k\geq 0} r^k \sum_{j=1}^{d_k} c_k^j Y^{(k)}_j(x')$ is harmonic in $\mathbb B$ and
satisfies the following condition
\begin{equation}\label{ngs}
N_s(g) = \sup_{0\leq\rho < 1} \sup_{y' \in \mathbb S} (1-\rho)^{\beta - \alpha + m + 1}
\left( \int_{\mathbb S} |\Lambda_{m+1}(g \ast P_{x'})(\rho y')|^s dx' \right)^{1/s} < \infty.
\end{equation}
\end{thm}

{\it Proof.} Since necessity of (\ref{ngs}) is contained in Lemma \ref{nec} we prove sufficiency of condition (\ref{ngs}). We assume $p$ and $q$ are finite, the remaining cases can be treated in a similar manner. Take $f \in B^{p, 1}_\alpha$ and set $h = M_c f$. Applying the operator $\Lambda_{m+1}$ to both sides of equation (\ref{simple}) we obtain
\begin{equation}\label{lamha}
\Lambda_{m+1} h(rx) = \int_{\mathbb S} \Lambda_{m+1} (g \ast P_{y'})(x) f(ry') dy'.
\end{equation}
Now we estimate the $L^s$ norm of the above function on $|x| = r$:
\begin{align}\label{form4}
M_s(\Lambda_{m+1} h, r^2) & \leq \int_{\mathbb S} M_s(\Lambda_{m+1} (g \ast P_{y'}), r) |f(ry')| dy' \notag \\
& \leq M_1(f, r) \sup_{y' \in \mathbb S} \left( \int_{\mathbb S} | \Lambda_{m+1} (g \ast P_{y'})(rx')|^s dx'
\right)^{1/s} \notag \\
& \leq M_1(f, r)N_s(g) (1-r)^{\alpha - \beta -m -1}.
\end{align}
Since,
$$\int_0^1 M_s^p(h, r^2)(1-r)^{\beta p -1} r^{n-1} dr  \leq C \int_0^1 (1-r)^{p(m+1)} M_s^p(\Lambda_{m+1}h, r^2) (1-r)^{\beta p -1}r^{n-1} dr,$$
see \cite{DS}, we have
\begin{align*}
\| h \|_{B^{p,s}_\beta}^p &  \leq C \int_0^1 (1-r)^{p(m+1)} M_s^p(\Lambda_{m+1}h, r^2) (1-r)^{\beta p -1}r^{n-1} dr \\
& \leq C N_s^p(g) \int_0^1 M_1^p(f, r)(1-r)^{\alpha p - 1} r^{n-1} dr \\
& = CN_s^p(g) \| f \|_{B^{p,1}_\alpha}^p,
\end{align*}
and therefore $\| h \|_{B^{p,s}_\beta} \leq \| f \|_{B^{p,1}_\alpha}$. Since $\| h \|_{B^{q,s}_\beta} \leq C
\| h \|_{B^{p,s}_\beta}$, see (\ref{inc}), the proof is complete. $\Box$

Next we develop another approach to multiplier problems in harmonic function spaces that hinges upon duality results. The first of these duality results is from \cite{Za1}: $(B^{p, q}_\alpha)^\ast \simeq B^{p', q'}_\alpha$, where
$1 < p < \infty$, $1 \leq q \leq \infty$ and $p'$ (resp. $q'$) is the exponent conjugate to $p$ (resp. $q$). Namely, the above identification of the dual space is with respect to the pairing
$$\langle f, g \rangle = \int_{\mathbb B} f(x)g(x) (1-|x|^2)^{2\alpha - 1} dx, \qquad f \in B^{p,q}_\alpha,
\quad g \in B^{p', q'}_\alpha.$$
In particular, the spaces $B^{p,q}_\alpha$ are reflexive for $1 < p < \infty$, $1 \leq q \leq \infty$. One can easily verify that for $c \in M_H(B^{p_1,q_1}_\alpha, B^{p_2, q_2}_\alpha)$, where $1 < p_1, p_2 < \infty$,
$1 \leq q_1, q_2 \leq \infty$, the adjoint operator $M_c^\ast : B^{p'_2, q'_2}_\alpha \rightarrow B^{p'_1, q'_1}_\alpha$ is also a multiplier operator generated by $c$.

\begin{thm}\label{bpinf}
Let $1 < p, q < \infty$ and $m > \alpha - 1$. Then for a double indexed sequence
$c = \{ c_k^j : k \geq 0, 1 \leq j \leq d_k \}$ the following conditions are equivalent:

1. $c \in M_H(B^{p,\infty}_\alpha, B^{q,\infty}_\alpha)$.

2. The function $g(x) = \sum_{k\geq 0} r^k \sum_{j=1}^{d_k} c_k^j Y^{(k)}_j(x')$ is harmonic in $\mathbb B$ and
satisfies the following condition
\begin{equation}\label{ng}
N_1(g)< \infty.
\end{equation}
\end{thm}

{\it Proof.} Assume $c \in M_H(B^{p,\infty}_\alpha, B^{q,\infty}_\alpha)$, then $M_c$ maps $B^{p,\infty}_\alpha$
continuously into $B^{q,\infty}_\alpha$ and therefore $M_c^\ast = M_c$ maps $B^{q',1}_\alpha$ into $B^{p',1}_\alpha$.
Hence, by Theorem \ref{bpbp}, condition (\ref{ng}) is satisfied. Conversely, assume (\ref{ng}) is satisfied. Then, by Theorem \ref{bpbp}, $M_c$ maps $B^{q',1}_\alpha$ into $B^{p',1}_\alpha$. Therefore, $M_c^\ast = M_c$ maps
$(B^{p,\infty}_\alpha)^{\ast\ast} = B^{p,\infty}_\alpha$ continuously into $(B^{q,\infty}_\alpha)^{\ast\ast} =
B^{q,\infty}_\alpha$, and the proof is complete. $\Box$

Now we recall a duality result from \cite{JP1}.
\begin{prop}
Let $1 < p < \infty$, $\alpha > 0$ and let $q$ be the exponent conjugate to $p$. Then $(A^p_\alpha)^\ast \simeq
A^q_\alpha$ with respect to the pairing
$$ \langle u, v \rangle = \int_{\mathbb B} u(x)v(x)(1-|x|^2)^\alpha dx.$$
\end{prop}

Again, it is easy to see that if we have a multiplier $M_c : A^{p_1}_\alpha \rightarrow A^{p_2}_\alpha$, where
$1 < p_1, p_2 < \infty$, then the adjoint operator $M_c^\ast : A^{p'_2}_\alpha \rightarrow A^{p'_1}_\alpha$ is
also a multiplier operator generated by the same sequence $c$. Therefore we have the following proposition.

\begin{prop}
$M_H(A^{p_1}_\alpha, A^{p_2}_\alpha) = M_H(A^{p'_2}_\alpha, A^{p'_1}_\alpha)$, $1 < p_1, p_2 < \infty$, $\alpha > 0$.
\end{prop}

The following duality result is also from \cite{JP1}.

\begin{prop}
For $\alpha > -1$ we have $B_0^\ast \simeq A^1_\alpha$ and $(A^1_\alpha)^\ast \simeq B$. Both identifications are
with respect to the following pairing:
\begin{equation}\label{pair}
\langle f, g \rangle = \int_{\mathbb B} f(x)g(x) (1-|x|^2)^\alpha dx.
\end{equation}
\end{prop}
Note that for given $g \in B$ integral in (\ref{pair}) is not necessarily convergent for all $f \in A^1_\alpha$, however, it converges on a dense subset $A^2_\alpha$ and extends by continuity to $A^1_\alpha$, see \cite{JP1} for
details. Since in this situation we again have $M_c^\ast = M_c$ one immediately obtains
\begin{equation}\label{subset}
M_H(B_0, B_0) \subset M_H(A^1_\alpha, A^1_\alpha) \subset M_H(B, B), \qquad \alpha > -1.
\end{equation}
This allows us to characterize multipliers from $B_0$ to $B_0$.

\begin{thm}
Let $m>-1$. For a double indexed sequence $c = \{ c_k^j : k \geq 0, 1 \leq j \leq d_k \}$ the following conditions are equivalent:

1. $c \in M_H(B_0, B_0)$.

2. The function $g(x) = \sum_{k\geq 0} r^k \sum_{j=1}^{d_k} c_k^j Y^{(k)}_j(x')$ is harmonic in $\mathbb B$ and
satisfies the following condition:
\begin{equation}\label{lbl}
\tilde N_1(g) = \sup_{0\leq\rho < 1} \sup_{y' \in \mathbb S} (1-\rho)^{m + 1}
\int_{\mathbb S} |\Lambda_{m+1}(g \ast P_{x'})(\rho y')| dx' < \infty.
\end{equation}
\end{thm}

{\it Proof.} Since $A^1_\alpha = B^{1,1}_{\alpha + 1}$, Theorem \ref{bpbp} and (\ref{subset}) show that condition
(\ref{lbl}) is necessary for $c \in M_H(B_0, B_0)$. Now we assume $g_c = g$ satisfies (\ref{lbl}) and choose
$f \in B_0$. Set $h = M_c f$. Since $\Lambda_{m+1} \nabla = \nabla \Lambda_{m+1}$ on $h(\mathbb B)$ we obtain, for $0 \leq r < 1$:

The proof is going to rely on the following easily checked identity
$$\int_{\mathbb S} \nabla^{-1}G(y) \nabla g(y) dy = \int_{\mathbb S} G(y)g(y) dy,$$
where $G: \mathbb B \rightarrow \mathbb C^n$ is a gradient of a harmonic function and $g \in h(\mathbb B)$.

Now we have
\begin{align*}
\Lambda_{m+1} \nabla h(rx) & = \nabla \Lambda_{m+1} h(rx) = \int_{\mathbb S} \nabla_x (\Lambda_{m+1}(g \ast P_{y'})
(x) f(ry') dy'\\
& = \int_{\mathbb S} \nabla_y^{-1} \nabla_x (\Lambda_{m+1}(g \ast P_{y'})(x) \nabla_y f(y) dy'.
\end{align*}

Therefore
\begin{align*}
|\Lambda_{m+1} \nabla h(rx)| & \leq (1-r)^{-1} \| f \|_{B_0} \int_{\mathbb S} |\Lambda_{m+1}
\nabla_y^{-1}\nabla_x (g\ast P_{y'})(rx')| dy'\\
& = (1-r)^{-1} \| f \|_{B_0} \int_{\mathbb S} |\Lambda_{m+1} \nabla_y^{-1} (g\ast \nabla_x P_{x'})(ry')| dy'\\
& = (1-r)^{-1} \| f \|_{B_0} \int_{\mathbb S} |\Lambda_{m+1} \nabla_y^{-1} (g\ast \nabla_y P_{x'})(ry')| dy'\\
& = (1-r)^{-1} \| f \|_{B_0} \int_{\mathbb S} |\Lambda_{m+1} (g\ast P_{x'})(ry')| dy'\\
& \leq \tilde N_1(g) \| f \|_{B_0}(1-r)^{-m-2}, \qquad r \in I, \quad x \in \mathbb B.
\end{align*}
Therefore $M_\infty(\Lambda_{m+1}\nabla h, r^2) \leq C (1-r)^{-m-2}$ for all $r \in I$. This implies, see \cite{DS},
Chapter 7, that $M_\infty(\nabla h, r ^2) \leq C (1-r)^{-1}$. This proves that $M_c$ maps $B_0$ into $B$. Since harmonic polynomials are dense in $B_0$, see \cite{JP2}, and $M_c$ maps harmonic polynomials into harmonic polynomials it follows that $M_c$ maps $B_0$ into $B_0$. $\Box$

The following proposition is a partial extension of Theorem \ref{bpinf}.
\begin{prop}\label{bpibpi}
Let $0<p\leq 1$, $p \leq q \leq \infty$ and $m > \alpha - 1$. If a double indexed sequence $c = \{ c_k^j : k \geq 0, 1 \leq j \leq d_k \}$ satisfies the following condition:
\begin{equation}\label{ng1}
N_1(g) = \sup_{0\leq\rho < 1} \sup_{y' \in \mathbb S} (1-\rho)^{m+1-\alpha + \beta} \int_{\mathbb S}
|\Lambda_{m+1}(g \ast P_{x'})(\rho y')| dx' < \infty,
\end{equation}
where $g(x) = \sum_{k\geq 0} r^k \sum_{j=1}^{d_k} c_k^j Y^{(k)}_j(x')$, then $c \in M_H(B^{p,\infty}_\alpha, B^{q,\infty}_\beta)$.
\end{prop}

{\it Proof.} Let us assume that condition (\ref{ng1}) is satisfied. Let $f \in B^{p,\infty}_\alpha (\mathbb B)$ and set $h = c \ast f$. We have, using Lemma \ref{intcon}:
\begin{align}
h(r\rho x') & = \int_{\mathbb S} (g \ast P_{x'})(ry') f(\rho y') dy' \notag  \\
& = 2 \int_0^1 \int_{\mathbb S} \Lambda_{m+1} (g \ast P_{x'})(rR \xi) f(\rho R \xi) (1-R^2)^m R^{n-1} d\xi dR.
\label{hrro}
\end{align}
Using (\ref{hrro}) and $M_\infty(f, \rho R) \leq M_\infty(f, R)$ we obtain, for $x = rx' \in \mathbb B$:
\begin{equation*}
|h(r\rho x')| \leq 2 \int_0^1 \int_{\mathbb S} |\Lambda_{m+1}(g \ast P_{x'})(rR\xi)| M_\infty(f, R) (1-R^2)^m
R^{n-1} d\xi dR.
\end{equation*}
Now letting $\rho \rightarrow 1$ and using condition (\ref{ng1}) we obtain
\begin{align*}
|h(rx')| & \leq 2 \int_0^1 M_\infty(f, R) (1-R^2)^m R^{n-1}\int_{\mathbb S} |\Lambda_{m+1}(g \ast P_{x'})(rR\xi)|
d\xi dR \\
& \leq 2^{m+1}  \int_0^1 M_\infty(f, R) (1-R)^m R^{n-1}\int_{\mathbb S} |\Lambda_{m+1}(g \ast P_{\xi})(rRx')|
d\xi dR \\
& \leq 2^{m+1} N_1(g) \int_0^1 M_\infty(f, R) (1-R)^m R^{n-1} (1-rR)^{\alpha - \beta - m -1} dR.
\end{align*}
Since $M_\infty(f, R)$ is an increasing function we can apply Lemma 3 from \cite{AS} to obtain
\begin{equation*}
|h(rx')|^p \leq C \int_0^1 M_\infty^p(f, R) \frac{(1-R)^{mp + p -1}}{(1-rR)^{p(m+1+\beta - \alpha)}} R^{n-1} dR.
\end{equation*}
This estimate is valid for all $x' \in \mathbb B$ and therefore gives estimate for $M_\infty^p(h, r)$ which is used,
together with Lemma \ref{rro}, in the following inequalities:
\begin{align*}
\| h \|_{B^{p, \infty}_\beta}^p & = \int_0^1 M_\infty^p(h, r)(1-r)^{\beta p -1} r^{n-1}dr \\
& \leq C \int_0^1 M_\infty^p(f, R) (1-R)^{mp+p-1}R^{n-1} \int_0^1
\frac{(1-r)^{\beta p -1}}{(1-rR)^{p(m+1+\beta - \alpha)}} r^{n-1}dr dR \\
& \leq C \int_0^1 M_\infty^p(f, R) (1-R)^{\alpha p -1} R^{n-1} dR = C \| f \|_{B^{p, \infty}_\alpha}^p.
\end{align*}
Since, by (\ref{inc}), $\| h \|_{B^{q, \infty}_\beta} \leq C \| h \|_{B^{p, \infty}_\beta}$ for $p \leq q$, the proof is complete. $\Box$

Our last theorem shows that restriction $1 \leq p$ in Theorem \ref{bpbp} can be removed in the case $s = \infty$.
\begin{thm}\label{bpbpi}
Let $0 < p \leq q \leq \infty$ and $m > \alpha - 1$. Then for a double indexed sequence
$c = \{ c_k^j : k \geq 0, 1 \leq j \leq d_k \}$ the following conditions are equivalent:

1. $c \in M_H(B^{p,1}_\alpha, B^{q,\infty}_\beta)$.

2. The function $g(x) = \sum_{k\geq 0} r^k \sum_{j=1}^{d_k} c_k^j Y^{(k)}_j(x')$ is harmonic in $\mathbb B$ and
satisfies the following condition
\begin{equation}\label{ngi}
N_\infty(g)< \infty.
\end{equation}
\end{thm}

{\it Proof.} The case $1 \leq p \leq \infty$ was settled in Theorem \ref{bpbp}. Assume $0<p<1$, the necessity of the
condition $N_\infty(g) < \infty$ was established in Lemma \ref{nec}. Now we assume $g = g_c$ satisfies $N_\infty(g) <
\infty$ and we use the same method of proof as in Proposition \ref{bpibpi}. Let $f \in B^{p,1}_\alpha$ and set
$h = M_c f$. Starting from (\ref{hrro}) and using $M_1(f, \rho R) \leq M_1(f, R)$ we obtain
$$|h(r\rho x')| \leq 2 \int_0^1 (1-R^2)^m R^{n-1} M_1(f, R) \sup_{\xi, x' \in \mathbb S} |\Lambda_{m+1} (g \ast P_{x'})
(rR \xi)| dR.$$
Then, analogously to the proof of Proposition \ref{bpibpi}, we let $\rho \rightarrow 1$ to obtain, using condition
(\ref{ngi}):
$$|h(rx')| \leq 2^{m+1} N_\infty(g) \int_0^1 M_1(f, R) (1-R)^m (1-rR)^{\alpha - \beta - m - 1} R^{n-1} dR.$$

Since $M_1(f, R)$ is an increasing function, one can follow the same reasoning as in the proof of Proposition
\ref{bpibpi}, replacing $M_\infty(f, R)$ by $M_1(f, R)$, to obtain $\| h \|_{B^{p, \infty}_\beta} \leq C
\| f \|_{B^{p,1}_\alpha}$. That suffices, due to embedding (\ref{inc}). $\Box$

\end{document}